\documentclass{article}

\usepackage{
 amsmath,
 amsxtra,
 amsthm,
 amssymb,
 etex,
 mathrsfs,
 stmaryrd,
pdflscape
 }
\usepackage[all]{xy}
\usepackage{hyperref}

\newtheorem{theorem}{Theorem}[section]
\newtheorem{lemma}[theorem]{Lemma}
\newtheorem{conjecture}[theorem]{Conjecture}
\newtheorem{proposition}[theorem]{Proposition}
\newtheorem{corollary}[theorem]{Corollary}

\theoremstyle{definition}
\newtheorem{defn}[theorem]{Definition}
\newtheorem{remark}[theorem]{Remark}

\newcommand{\bd}{\begin{defn}}
\newcommand{\ed}{\end{defn}}
\newcommand{\bl}{\begin{lemma}}
\newcommand{\el}{\end{lemma}}
\newcommand{\bp}{\begin{proposition}}
\newcommand{\ep}{\end{proposition}}
\newcommand{\bt}{\begin{theorem}}
\newcommand{\et}{\end{theorem}}
\newcommand{\bc}{\begin{corollary}}
\newcommand{\ec}{\end{corollary}}
\newcommand{\br}{\begin{remark}}
\newcommand{\er}{\end{remark}}
\newcommand{\ba}{\begin{array}}
\newcommand{\ea}{\end{array}}
\newcommand{\bpf}{\begin{proof}}
\newcommand{\epf}{\end{proof}}

\newcommand{\Z}{\mathbb{Z}}
\newcommand{\Q}{\mathbb{Q}}
\newcommand{\Zp}{\mathbb{Z}_{p}}
\newcommand{\Qp}{\mathbb{Q}_{p}}

\newcommand{\Op}{\mathcal{O}}

\newcommand{\al}{\alpha}
\newcommand{\be}{\beta}
\newcommand{\Ga}{\Gamma}

\newcommand{\La}{\Lambda}
\newcommand{\la}{\lambda}

\newcommand{\Iw}{\mathrm{Iw}}
\newcommand{\p}{\mathfrak{p}}

\newcommand{\Sel}{\mathrm{Sel}_{p^\infty}} \DeclareMathOperator{\Gal}{Gal}
 \DeclareMathOperator{\rank}{rank}
\DeclareMathOperator{\corank}{corank}
\DeclareMathOperator{\Ext}{Ext}

\newcommand{\ord}{\mathrm{ord}}
\newcommand{\cyc}{\mathrm{cyc}}
\newcommand{\ac}{\mathrm{ac}}
\newcommand{\ch}{\mathrm{char}}

\newcommand{\M}{\mathfrak{M}}
\newcommand{\N}{\mathfrak{N}}

\newcommand{\ot}{\otimes}
\newcommand{\ilim}{\displaystyle \mathop{\varinjlim}\limits}
\newcommand{\plim}{\displaystyle \mathop{\varprojlim}\limits}

\newcommand{\coker}{\mathrm{coker}\,}

\newcommand{\lra}{\longrightarrow}

\newcommand{\ps}[1]{\llbracket #1 \rrbracket}

  \DeclareFontFamily{U}{wncy}{}
  \DeclareFontShape{U}{wncy}{m}{n}{<->wncyr10}{}
  \DeclareSymbolFont{mcy}{U}{wncy}{m}{n}
  \DeclareMathSymbol{\sha}{\mathord}{mcy}{"58}

\begin{document}
\title{On Iwasawa theory of abelian varieties over $\Zp^2$-extension with applications to Diophantine stability and integally Diophantine extensions}
 \author{
   Meng Fai Lim\footnote{School of Mathematics and Statistics, Key Laboratory of Nonlinear Analysis and Applications (Ministry of Education),
Central China Normal University, Wuhan, 430079, P.R.China.
 E-mail: \texttt{limmf@ccnu.edu.cn}}
}
\date{}
\maketitle

\begin{abstract} \footnotesize
\noindent We present certain results on the Iwasawa theory of an abelian variety with potentially good ordinary reduction at all primes above $p$. These are then applied to study
Diophantine stability and integally Diophantine extensions. Along the way, we also obtain some results pertaining to Mazur growth conjecture which refine previous results of Gajek-Leonard, Hatley, Kundu and Lei. Finally, we extend our investigation to the case of an elliptic curve with good supersingular reduction at the prime $p$ and make a similar analysis.


\end{abstract}

\section{Introduction}

Throughout the paper, $p$ will always denote a fixed odd prime. Let $A$ be an abelian variety defined over a number field $F$ which has potentially good ordinary reduction at all primes of $F$ above $p$. A well-known conjecture of Mazur \cite{Maz} predicts that the ($p$-primary) Selmer group of the said abelian variety over the cyclotomic $\Zp$-extension $F_\cyc$ of $F$ is cotorsion over the Iwasawa algebra of $\Ga=\Gal(F_\cyc/F)$. The works of Rubin \cite{Ru} and Kato \cite{K} provides many evidences towards this conjecture.

The aim of this paper is study certain consequences of this conjecture of Mazur. Denote by $X(A/F_\cyc)$ the Pontryagin dual of the Selmer group of $A$ over $F_\cyc$. Assuming Mazur conjecture conjecture, we may appeal to the structure theory of finitely generated torsion $\Zp\ps{\Ga}$-module and Weierstrass preparation Theorem to define the characteristic polynomial of $X(A/F_\cyc)$, which we write as $\ch_{\Zp\ps{\Ga}}X(A/F_\cyc)$. We shall identity $\Zp\ps{\Ga}$ with the power series ring $\Zp\ps{U}$, and under this identification, it makes sense to speak of $\ord_U\big(\ch_{\Zp\ps{\Ga}}X(A/F_\cyc)\big)$ which is the highest power of $U$ dividing $\ch_{\Zp\ps{\Ga}}X(A/F_\cyc)$. We can now present our first main result.

\bt[Theorem \ref{main thm Sel}]
Suppose that $A$ is an abelian variety defined over $F$ which has potential good ordinary reduction at every prime of $F$ above $p$. Let $F_\infty$ be a $\Zp^2$-extension of $F$ which contains the cyclotomic $\Zp$-extension $F_\cyc$. Suppose that $X(A/F_\cyc)$ is torsion over $\Zp\ps{\Gal(F_\cyc/F)}$. Denote by $\mathfrak{M}(A/F_\infty)$ the set of $\Zp$-extensions $K_\infty$ of $F$ contained in $F_\infty$ such that the Mordell-Weil rank of $A$ is unbounded in $K_\infty$. Then the set $\mathfrak{M}(A/F_\infty)$ is finite with
\[ \#\mathfrak{M}(A/F_\infty) \leq \ord_U\left(\ch_{\Zp\ps{\Ga}}\big(X(A/F_\cyc)\big)\right).
\]\et

When $A$ is an elliptic curve, the above result has been established in \cite[Proposition 3.9]{GHKL}, albeit under an extra technical assumption. Our result here not only generalizes this to an abelian variety but also remove the said extra assumption (see Remark \ref{main thm Sel remark}(1)). As a consequence, we may apply our theorem to refine some of the results of Gajek-Leonard-Hatley-Kundu-Lei \cite{GHKL} and Kundu-Lei \cite{KL} pertaining to the Mazur growth conjecture.

We next come to the study of Diophantine stability and integally Diophantine extensions. Here the notion of Diophantine stability is defined as in the sense of Mazur-Rubin \cite{MRL}. On the other hand, the integrally Diophantine extensions are intimately related with the Hilbert's 10th problem (see \cite{MR}; also see \cite{GP, KP, KLS, Ra} and references therein for further works on this).
We will apply our theorem to show that, under the assumption of Mazur's conjecture, there exist an abundance of cyclic extensions, where the abelian variety is Diophantine stable. Consequently, we then see that there is an abundance of cyclic extensions, which are integally Diophantine extensions. However, we must emphasize that while our results show the existence of many such extensions, whether pertaining to Diophantine stability or integrally Diophantine extensions, we are currently unable to ensure that they are specifically extensions of the base field $F$.
To provide a more precise description of our result in these aspects, we begin by letting $n$ denote a given positive integer. Denote by $\mathfrak{D}_n(A/F)$ the set consisting of pairs $(L_1, L_2)$ satisfying the following:
 \begin{enumerate}
   \item[$(1)$] $L_1$ and $L_2$ are contained in a $\Zp$-extension of $F$ with $|L_1:L_2|=p^n$.
   \item[$(2)$] The abelian variety $A$ is diophantine-stable for $L_1/L_2$.
 \end{enumerate}

We also denote by $\mathcal{I}_n(F)$ the set which consists of pairs $(L_1, L_2)$ satisfying the following:

 \begin{enumerate}
   \item[$(1)$] $L_1$ and $L_2$ are contained in a $\Zp$-extension of $F$ with $|L_1:L_2|=p^n$.
   \item[$(2)$] $L_1/L_2$ is integrally Diophantine.
 \end{enumerate}

Our results are then as follows.

\bt[Theorem \ref{DioStable}]
Let $F$ be a number field with at least one complex prime.
Suppose that $A$ is an abelian variety defined over $F$ with potential good ordinary reduction at every prime above $p$. Assume that $X(A/F_\cyc)$ is torsion over $\Zp\ps{\Gal(F_\cyc/F)}$. Then the set
$\mathfrak{D}_n(A/F)$ is uncountable for every integer $n\geq 1$.
\et

\bt[Theorem \ref{intDio}]
Let $F$ be a number field with at least one complex prime.
Suppose that there exist an abelian variety $A$ defined over $F$ which satisfies all of the following.
\begin{enumerate}
  \item[$(i)$] The abelian variety $A$ has potential good ordinary reduction at every prime above $p$.
  \item[$(ii)$] $X(A/F_\cyc)$ is torsion over $\Zp\ps{\Gal(F_\cyc/F)}$.
  \item[$(iii)$] $\rank_\Z A(F) >0$.
\end{enumerate}
Then, for every integer $n\geq 1$, the set
$\mathcal{I}_n(F)$ is uncountable.
\et

Combining our result with the deep works of Rubin and Kato, we therefore have the following unconditional result on the Diophantine stability of an elliptic curve defined over $\Q$.

\bc[Corollary \ref{DioStable coro}]
Let $E$ be an elliptic curve defined over $\Q$ with good ordinary reduction at $p$.
Suppose that $F$ is a finite abelian imaginary extension of $\Q$. Then for every integer $n\geq 1$, the set
$\mathfrak{D}_n(E/F)$ is uncountable.
\ec

Furthermore, by invoking deep results of Rubin, Kato, Bump-Friedberg-Hoffstein, Murty-Murty and Gross-Zagier, we have the following unconditional conclusion on the set $\mathcal{I}_n(F)$.

\bc[Corollary \ref{intDio coro}]
Suppose that $F$ is a finite abelian imaginary extension of $\Q$. Then for every integer $n\geq 1$, the set
$\mathcal{I}_n(F)$ is uncountable.
\ec

Recall that a conjecture of Denef and Lipshitz predicts that $K/\Q$ is integrally diophantine for every number field $K$ (see \cite{DL}). Recently, Koymans and Pagano announced a proof of this conjecture (see \cite{KP}). Taking our Theorem \ref{intDio} and Corollary \ref{intDio coro} into account, one is tempted to formulate the following conjecture.

\medskip \noindent \textbf{Relative Denef-Lipshitz conjecture.}
Every finite extension of number field $L/K$ is integrally diophantine.

\medskip
We also note that in the works \cite{GP, KLS}, the authors actually established new cases of Denef-Lipshitz conjecture by proving appropriate cases of the relative Denef-Lipshitz conjecture.

So far, the discussion of the paper revolves around the situation when the abelian variety in question has potential good ordinary reduction at every primes above $p$. It's a natural question to ask what can be said if the abelian variety does not have good ordinary reduction at primes above $p$. In response to this, we present some results in the following modest context.

Hereafter, $E$ will denote an elliptic curve over $\Q$ which has good supersingular reduction at the prime $p$. Denote by $\widetilde{E}$ the reduced curve of $E$ modulo $p$. We will always assume that $a_p= p+1 - |\widetilde{E}(\mathbb{F}_p)|=0$ (note that this automatically holds if $p\geq 5$). Let $F$ be an imaginary quadratic field of $\Q$ at which the prime $p$ splits completely, say $p=\p\overline{\p}$. Let $F_\infty$ be the $\Zp^2$-extension of $F$. Following the ideas of Kobayashi \cite{Kob}, B. D. Kim defined his multi-signed Selmer groups in his paper \cite{Kim}, which we will address them as ``signed Selmer groups'' and denote their Pontryagin dual by $X^{s,z}(E/F_\infty)$ for $s,z\in\{+,-\}$ (see body of paper for their precise definitions).

In this context, the following is the analog of Theorem \ref{main thm Sel}.

\bt[Theorem \ref{signed X leq ord}]
Let $E$ be an elliptic curve of conductor $N$ over $\Q$ with good supersingular reduction at the prime $p$ and $a_p=0$. Let $F$ be an imaginary quadratic field of $\Q$ at which all the prime divisors of $p$ split completely in $F/\Q$. If $X^{s,z}(E/F_\cyc)$ is torsion over $\Zp\ps{\Ga}$ for every $s,z\in\{+,-\}$, then we have
\[ \#\mathfrak{M}(E/F_\infty) \leq
 \sum_{s,z\in\{+,-\}}\ord_U\big(\ch_{\Zp\ps{\Ga}}X^{s,z}(E/F_\cyc)\big).
\]
\et

Unfortunately, at present, it's only known that $X^{+,+}(E/F_\cyc)$ and $X^{-,-}(E/F_\cyc)$   are torsion over $\Zp\ps{\Ga}$ (see \cite{Kob, LeiP}). Nevertheless, we can at least have the following unconditional result concerning Diophantine stability.

\bt[Theorem \ref{supersingular thm}]
Let $E$ be an elliptic curve of conductor $N$ over $\Q$ with good supersingular reduction at the prime $p$ and $a_p=0$. Let $F$ be an imaginary quadratic field of $\Q$ such that  $p$ splits completely in $F/\Q$, and the primes of $F$ above $p$ are totally ramified in $F_{\ac}/F$. Suppose either of the following statement is valid.
\begin{enumerate}
  \item[$(i)$] $\Sel(E/F)$ is finite.
  \item[$(ii)$] Every prime divisor of $N$ splits completely in $F/\Q$.
  \end{enumerate}
Then, for every integer $n\geq 1$, the set
$\mathfrak{D}_n(E/F)$ is uncountable.
\et

We now provide a brief overview of the paper. In Section \ref{Algebraic section}, we collect certain algebraic results which will be required for our arithmetic discussion. Moving on to Section \ref{main thm section}, we will prove our main result (Theorem \ref{main thm Sel}) on bounding the number of $\Zp$-extensions at which the Mordell-Weil rank is unbounded. Subsequently, in Section \ref{Mazur Growth section}, we apply this said theorem to improve the results of \cite{GHKL,KL} pertaining to Mazur Growth Conjecture. Section \ref{DioStable section} will be devoted to exploring the consequences of our main theorem on Diophantine stability and integrally Diophantine extensions. We then make some remark on Theorem \ref{main thm Sel} in Section \ref{Remark on Thm3.3 section}. Finally, in Section \ref{elliptic supersingular section}, we establish results analogue to those in Sections \ref{main thm section} and \ref{DioStable section} for an elliptic curve with good supersingular reduction over the $\Zp^2$-extension of an imaginary quadratic field.

\subsection*{Acknowledgments}
We like to thank Antonio Lei for his interest and comments on the paper.
The author is partially supported by the Fundamental Research Funds for the Central Universities No.\ CCNU25JCPT031, and the Open Research Fund of Hubei Key Laboratory of Mathematical Sciences (Central China Normal University).

\section{Algebraic preliminaries} \label{Algebraic section}

In this section, we establish certain algebraic preliminaries
and notation that are necessary for the discussion in the paper. We always let $G$ denote the group $\Zp^2$, and denote by $\La$ the Iwasawa algebra $\Zp\ps{G}$. We fix a subgroup $H$ of $G$ such that $H\cong\Zp\cong G/H$. Fix a choice of topological generators $\sigma$ and $\tau$ of $G$ such that $\tau$ is a topological generator of $H$ and $\sigma~(\mathrm{mod}~H)$ is a topological generator of $\Ga:=G/H$. We then identified $\Zp\ps{G}$ with the power series ring $\Zp\ps{U,W}$ under the correspondence $\sigma\mapsto U, \tau\mapsto W$. We will also frequently identify $\Zp\ps{\Ga}$ with the power series ring $\Zp\ps{U}$.

For each $(a,b)\in \mathbb{P}_1(\Zp)$, we denote by $H_{a,b}$ the subgroup of $G$ topologically generated by
$\sigma^a\tau^b$. In particular, one has $H= H_{0,1}$. The Iwasawa algebra $\Zp\ps{H_{a,b}}$ can then be identified as the power series ring $\Zp\ps{(U+1)^a(W+1)^b-1}$. Moreover, we have $G/H_{a,b}\cong \Zp$. We then write $\Ga_{a,b} = G/H_{a,b}$, and in particularly, write $\Ga = \Ga_{0,1}=G/H$.

\bd
Let $M$ be a torsion $\Zp\ps{G}$-module. Set
\[\mathfrak{A}(M): = \{(a,b)\in P_1(\Zp) \mid M_{H_{a,b}}~ \mbox{is not torsion over}~\Zp\ps{G/H_{a,b}}\}.\]
\ed

A commutative algebraic argument tells us that the set $\mathfrak{A}(M)$ is finite (for instance, see \cite[Lemma 2.2]{LimFineMod}). The aim of this section is to give an upper bound of the cardinality of this set which is the content of the following proposition.

\bp \label{N leq ord}
Let $M$ be a torsion $\Zp\ps{G}$-module with the property that $M_H$ is torsion over $\Zp\ps{\Ga}$. Then we have
\[ \#\mathfrak{A}(M) \leq \ord_{U}(\ch_{\Zp\ps{\Ga}}M_H). \]
\ep

Note that $M_H = M/W$ and $H_1(H,M) = M[W]$, where the latter is the submodule of $M$ consisting of elements annihilated by $W$. We will make use of these identifications without further mention. The remainder of the section will be devoted to the proof of the preceding proposition. As a start, we record the following lemma which is taken from a monograph of Perrin-Riou \cite{PR}.

\bl \label{charproject}
Let $M$ be a torsion $\Zp\ps{G}$-module with the property that $M_H$ is torsion over $\Zp\ps{\Ga}$. Denote by $\pi$ the natural projection map $\Zp\ps{G}\twoheadrightarrow\Zp\ps{\Ga}$. Then $H_1(H,M)$ is also torsion over $\Zp\ps{\Ga}$, and we have
\[ \ch_{\Zp\ps{\Ga}}(H_1(H,M)) \cdot \pi\big(\ch_{\Zp\ps{G}}(M)\big) = \ch_{\Zp\ps{\Ga}}(M_H). \]
In particular, if $H_1(H,M)$ is pseudo-null as a $\Zp\ps{\Ga}$-module, we then have
\[ \pi\big(\ch_{\Zp\ps{G}}(M)\big) = \ch_{\Zp\ps{\Ga}}(M_H). \]
\el

\bpf See \cite[\S I.1.3, Lemma 4]{PR}. %Alternatively, this is a special case of \cite[Lemma 2.8]{HaoL}. Nevertheless, we provide a direct proof here for the convenience of our readers.  Denote by $M_{\mathrm{null}}$ the maximal pseudo-null $\Op\ps{G}$-submodule of $M$. We then have the following exact sequence\[ 0 \lra M_{\mathrm{null}}\lra M \lra E \lra N\lra 0, \]where $N$ is some pseudo-null $\Op\ps{G}$-module and $E = \bigoplus_{\mathfrak{p}}\Op\ps{G}/\mathfrak{p}^{a_\mathfrak{p}}$ for some height one prime $\mathfrak{p}$ of $\Op\ps{G}$ not equal to $(W)$ (this is necessary due to the $\Op\ps{\Ga}$-torsionness of $M_H$). Denoting by $Z$ the image of the map $M\lra E$, we then have two short exact sequences\[ 0\lra M_{\mathrm{null}}\lra M  \lra Z\lra 0, \]
%\[ 0\lra Z\lra E \lra N\lra 0.\]
%Since the prime ideals appearing in $E$ are not equal to $(W)$, we have $E[W]= H_1(H,E)=0$ which in turn implies that $H_1(H,Z)=0$. Taking these into account, we obtain the following two exact sequences
%\[ 0\lra (M_{\mathrm{null}})_H\lra M_H  \lra Z_H\lra 0, \]
%\[ 0\lra N[W] \lra Z/W\lra E/W \lra N/W \lra 0.\]
%Now since $N$ is a pseudo-null $\Op\ps{G}$-module, it follows from \cite[Lemma 3.1]{Och} that $\ch_{\Op\ps{\Ga}}(N[W])= \ch_{\Op\ps{\Ga}}(N/W)$. Thus, it follows from this and the above exact sequence that we have $\ch_{\Op\ps{\Ga}}(Z/W)= \ch_{\Op\ps{\Ga}}(E/W)$. On the other hand, since $M[W]=0$ by assumption, we also have $M_{\mathrm{null}}[W]=0$. By \cite[Lemma 3.1]{Och}, this in turn implies that $M_{\mathrm{null}}/W$ is pseudo-null over $\Op\ps{\Ga}$. Therefore, the above short exact sequence yields $\ch_{\Op\ps{\Ga}}(M/W)= \ch_{\Op\ps{\Ga}}(Z/W)$, and hence we have the identity
%\[ \ch_{\Op\ps{\Ga}}(M/W)= \ch_{\Op\ps{\Ga}}(E/W) = \pi\big( \ch_{\Op\ps{G}}(E)\big).\]
%This is precisely what we set out to show.
\epf

We can now give the proof of Proposition \ref{N leq ord}.

\bpf[Proof of Proposition \ref{N leq ord}]
By definition, we have $(a,b)\in \mathfrak{A}(M)$ if and only if $M_{H_{a,b}}$ is not torsion over $\Zp\ps{\Ga_{a,b}}$. The latter will force the prime element $f_{a,b}:=(U+1)^a(W+1)^b-1$ to
  lie in the support of $M$ which in turn implies that $f_{a,b}$ divides $\ch_{\Zp\ps{G}}(M)$. Since the $f_{a,b}$'s are mutually coprime and $\Zp\ps{G}$ is a unique factorization domain, we have
  \[ \prod_{(a,b)\in\N(M)} f_{a,b}~\Big|~ \ch_{\Zp\ps{G}}(M). \]
  Applying the map $\pi$, we obtain the divisibility
  \[ \prod_{(a,b)\in\N(M)}(U+1)^a-1 ~\Big|~ \pi\big(\ch_{\Zp\ps{G}}(M)\big). \]
  Since $U~\Big|~(U+1)^a-1$, and taking Lemma \ref{charproject} into account, we obtain
  \[ U^{\#\N(M)} ~\Big|~ \ch_{\Zp\ps{\Ga}}(M_H). \]
  The conclusion of the proposition is now immediate from this.
  \epf

\section{Main theorem} \label{main thm section}

Let $A$ be an abelian variety defined over $F$ which has potential good ordinary reduction at every prime of $F$ above $p$. Let $S$ be a finite set of primes of $F$ which contains all the primes above $p$, the infinite primes and the primes of bad reduction of $A$. Denote by $F_S$ the maximal algebraic extension of $F$ which is unramified outside S. If $\mathcal{L}$ is a (possibly infinite) extension of $F$ contained in $F_S$, we write $G_S(\mathcal{L})=\Gal(F_S/\mathcal{L})$ and denote by $S(\mathcal{L})$ the set of primes of $\mathcal{L}$ above $S$.

For each $v\in S$ and a finite extension $L$ of $F$, set
\[ J_v(A/L)=\bigoplus_{w|v}H^1(L_w, A)_{p^\infty}.\]
If $\mathcal{L}$ is an infinite extension of $F$ contained in $F_S$, we define
\[ J_v(A/\mathcal{L})=\ilim_L J_v(A/L),\]
where $L$ runs through all the finite extensions of $F$ contained in $\mathcal{L}$.

The classical ($p$-primary) Selmer group of $A$ over $\mathcal{L}$ is defined by
\[ \Sel(A/\mathcal{L}) = \ker\left(H^1(G_S(\mathcal{L}), A_{p^\infty})\lra \bigoplus_{v\in S}J_v(A/\mathcal{L})\right).\]
 The Pontryagin dual of $\Sel(A/\mathcal{L})$ is then denoted
by $X(A/\mathcal{L})$.

Throughout our discussion, we shall frequently call upon the following conjecture of Mazur \cite{Maz}.

\medskip \noindent \textbf{Conjecture (Mazur).} Suppose that $A$ is an abelian variety defined over $F$ which has potential good ordinary reduction at every prime of $F$ above $p$. Let $F_\cyc$ be the cyclotomic $\Zp$-extension of $F$. Then $X(A/F_{\cyc})$ is a torsion $\Zp\ps{\Gal(F_\cyc/F)}$-module.

\medskip
Deep results of Kato \cite[Theorem 17.4]{K} and Rubin \cite[Theorem 4.4]{Ru} showed that the conjecture holds, whenever $A$ is an elliptic curve defined over $\Q$ with good ordinary reduction and $F$ is an abelian extension of $\Q$. Prior to these, Mazur established the validity of his conjecture under the assumption that $\Sel(A/F)$ is finite. In his proof, Mazur made use of the following, which nowadays is coined as ``Mazur control theorem''.

\bt[Mazur Control Theorem] \label{Mazur Control}
Let $A$ be an abelian variety defined over $F$ which has potential good ordinary reduction at every prime of $F$ above $p$. Suppose that $\mathcal{F}$ is a $\Zp$-extension of $F$ and $\mathcal{F}_n$ the intermediate subfield of $\mathcal{F}$ with $|\mathcal{F}_n:F|=p^n$. Then the restriction map
\[ \Sel(A/\mathcal{F}_n) \lra \Sel(A/\mathcal{F})^{\Gal(\mathcal{F}/\mathcal{F}_n)} \]
has finite kernel and cokernel which are bounded independent of $n$.
\et

\bpf
See \cite[Proposition 3.7]{G99}.
\epf

\bc \label{Mazur Control cor}
Let $A$ be an abelian variety defined over $F$ which has potential good ordinary reduction at every prime of $F$ above $p$. Let $\mathcal{F}$ be a $\Zp$-extension of $F$ and $\mathcal{F}_n$ the intermediate subfield of $\mathcal{F}$ with $|\mathcal{F}_n:F|=p^n$. Suppose that $X(A/\mathcal{F})$ is torsion over $\Zp\ps{\Gal(\mathcal{F}/F)}$. Then $n\to\infty$, we have
\[ \rank_\Z A(\mathcal{F}_n) = O(1). \]
\ec

\bpf
Since $X(A/\mathcal{F})$ is torsion over $\Zp\ps{\Gal(\mathcal{F}/F)}$, we have
\[ \mathrm{corank}_{\Zp} \left(\Sel(A/\mathcal{F})^{\Gal(\mathcal{F}/\mathcal{F}_n)} \right)= O(1). \]
By Theorem \ref{Mazur Control}, this in turn implies that
\[ \mathrm{corank}_{\Zp} \big(\Sel(A/\mathcal{F}_n)\big) = O(1). \]
Since we always have $\rank_\Z \big(A(\mathcal{F}_n)\big)\leq \mathrm{corank}_{\Zp} \big(\Sel(A/\mathcal{F}_n)\big)$, the conclusion of the corollary follows.
\epf

Let $F_\infty$ be a $\Zp^2$-extension of $F$ which contains the cyclotomic $\Zp$-extension $F_\cyc$.
Write $G=\Gal(F_\infty/F)\cong \Zp^2$. We identify the ring $\Zp\ps{G}$ with the power series ring $\Zp\ps{U,W}$ such that $\Zp\ps{\Gal(F_\infty/F_\cyc)}\cong \Zp\ps{W}$ and $\Zp\ps{\Gal(F_\cyc/F)}\cong \Zp\ps{U}$. Under this identification, we write $H_{a,b}$ for the subgroup of $G$ generated by $(1+U)^a(1+W)^b-1$, and $\Ga_{a,b}=G/H_{a,b}$. Denote by $F_{a,b}$ the fixed field of $H_{a,b}$. We are interested in the following set
\[\M(A/F_\infty) = \big\{(a,b)\in\mathbb{P}_1(\Zp)~|~\mbox{the Mordell-Weil rank of $A$ in $F_{a,b}$ is unbounded}\big\} \]

The main theorem of the section is presented below, where we show that the set is finite and provides an upper bound for the size of this particular set.

\bt \label{main thm Sel}
Let $A$ be an abelian variety defined over $F$ which has potential good ordinary reduction at every prime of $F$ above $p$. Let $F_\infty$ be a $\Zp^2$-extension of $F$ which contains the cyclotomic $\Zp$-extension $F_\cyc$. Suppose that $X(A/F_\cyc)$ is torsion over $\Zp\ps{\Gal(F_\cyc/F)}$. Then we have
\[ \#\mathfrak{M}(A/F_\infty) \leq \ord_U\left(\ch_{\Zp\ps{\Ga}}\big(X(A/F_\cyc)\big)\right).
\]\et

The remainder of the section will be devoted to the proof of Theorem \ref{main thm Sel}. As a start, we have the following general observation.

\bl \label{control all}
For every $(a,b)\in \mathbb{P}_1(\Zp)$,  the restriction map
\[ \Sel(A/F_{a,b}) \lra \Sel(A/F_\infty)^{H_{a,b}}\]
has kernel which is cofinitely generated over $\Zp$.
\el

\bpf
Indeed, the kernel of the restriction map is contained in $H^1(H_{a,b}, A_{p^\infty}(F_\infty))$, and the latter is plainly cofinitely generated over $\Zp$.
\epf

In the context of $F_\cyc$, we have the following sharper result on the restriction map.

\bl \label{control theorem}
The restriction map
\[ \Sel(A/F_\cyc) \lra \Sel(A/F_\infty)^H\]
has finite kernel and cokernel.
  \el

\bpf
Consider the following commutative diagram
 \[   \entrymodifiers={!! <0pt, .8ex>+} \SelectTips{eu}{}\xymatrix{
    0 \ar[r] &\Sel(A/F_\cyc)  \ar[d]^{\al} \ar[r] &  H^1\big(G_{S}(F_\cyc),A_{p^\infty}\big)
    \ar[d]^{h} \ar[r] & \displaystyle \bigoplus_{v\in S}J_v(A/F_\cyc) \ar[d]^{\oplus g_v}  \\
    0 \ar[r]^{} & \Sel(A/F_\infty)^{H} \ar[r]^{} & H^1\big(G_{S}(F_\infty),A_{p^\infty}\big)^{H} \ar[r] &
    \left(\displaystyle \bigoplus_{v\in S}J_v(A/F_\infty)\right)^{H}   } \]
    with exact rows. Via the snake lemma, it suffices to show that $\ker h$, $\ker g_v$ and $\coker h$ are finite. To begin with, we show that $h$ is surjective with a finite kernel. Indeed, since $H\cong \Zp$, the
restriction-inflation sequence tells us that the map $h$ is surjective with kernel $H^1(H, A_{p^\infty}(F_\infty))$. On the other hand, we have
\[ 0 = \rank_{\Zp\ps{H}}\Big( A_{p^\infty}(F_{\infty})\Big)^\vee = \rank_{\Zp} \Big(\big(A_{p^\infty}(F_{\infty})\big)^\vee\Big)_H - \rank_{\Zp} \Big(A_{p^\infty}(F_{\infty})^\vee\Big)^H\]
where the second equality follows from \cite[Proposition 5.3.20]{NSW}. But observe that
\[\Big(A_{p^\infty}(F_{\infty})^\vee\Big)_H = \Big(A_{p^\infty}(F_{\infty})^H\Big)^\vee = \Big(A_{p^\infty}(F_\cyc)\Big)^\vee,  \]
and the latter is finite by a theorem of Imai (\cite{Im}; also see \cite[Theorem 3.4]{Win}). Hence $\Big(A_{p^\infty}(F_{\infty})^\vee\Big)^H$ is finite. But this group is precisely
$H^1(H, A_{p^\infty}(F_\infty))^\vee$, and so we have our claim.

It remains to show that $g_v$ has finite kernel for every $v$. As a start, recall that for  $\mathcal{L}\in\{F_\cyc, F_\infty\}$, we have isomorphisms
\[ J_v(E/\mathcal{L})\cong \begin{cases} \ilim_{\mathcal{L}'}\bigoplus_{w|v}H^1(\mathcal{L}'_w, \tilde{A}_{w,p^\infty}),  & \mbox{if $v$ divides $p$}, \\
\ilim_{\mathcal{L}'}\bigoplus_{w|v}H^1(\mathcal{L}'_w, A_{p^\infty}), & \mbox{if $v$ does not divide $p$},\end{cases} \]
where the direct limit is taken over all finite extensions $\mathcal{L}'$ of $F_\cyc$ contained in $\mathcal{L}$ (see \cite[Propositions 4.1, 4.7 and 4.8]{CG}).

Taking these isomorphisms into account, it follows from the restriction-inflation sequence that we have $\ker g_v = \oplus_{w|v} H^1(H_w, D(F_{\infty,w}))$, where here the sum runs over all the primes of $F_c$ above $v$, $H_w$ is the decomposition group of $w$ in $H$ and $D$ denotes either $\tilde{A}_{w,p^\infty}$ or $A_{p^\infty}$ according as $v$ divides $p$ or not.
In particular, if $H_w=1$, then $H^1(H_w, D(F_{\infty,w}))=0$. This is indeed the case when $v$ does not divide $p$. It remains to consider primes $w$ which divides $p$ and that $H_w$ is nontrivial. Since $H\cong\Zp$, it then follows that $H_w\cong \Zp$. We may then apply the same argument as in the preceding paragraph to conclude that $H^1(H_w, \tilde{A}_{w,p^\infty})$ is finite. The proof of the proposition is now complete.
\epf

We are in position to prove Theorem \ref{main thm Sel}.

\bpf[Proof of Theorem \ref{main thm Sel}]
 Taking the hypothesis that $X(A/F_\cyc)$ is torsion over $\Zp\ps{\Ga}$ into account, it then follows from Lemma \ref{control theorem} that $X(A/F_\infty)_H$ is torsion over $\Zp\ps{\Ga}$
with
\[\ch_{\Zp\ps{\Ga}}X(A/F_\cyc) = \ch_{\Zp\ps{\Ga}}X(A/F_\infty)_H. \]
By \cite[Lemma 2.6]{HO}, $X(A/F_\infty)$ is a torsion $\Zp\ps{G}$-module. Thus, we may invoke Proposition \ref{N leq ord} to conclude that
\[ \#\mathfrak{A}\big(X(A/F_\infty)\big)\leq \ord_U\big(\ch_{\Zp\ps{\Ga}}X(A/F_\cyc)\big).\]
On the other hand, taking Corollary \ref{Mazur Control cor} and Lemma \ref{control all} into account, we have the following inclusion
\[ \mathfrak{M}\big(A/F_\infty\big) \subseteq \mathfrak{A}\big(X(A/F_\infty)\big). \]
The conclusion of theorem then follows from combining this inclusion with the above inequality.
\epf

\br \label{main thm Sel remark}
\begin{enumerate}
  \item[$(1)$] Our Theorem \ref{main thm Sel} strengthens \cite[Proposition 3.9]{GHKL}, where we remove the restrictive assumption that $X(A/F_\infty)$ is a direct sum of cyclic torsion modules.
  \item[$(2)$] Now, if $F_{a,b}$ is some $\Zp$-extension of $F$ with the properties that $A(F_{a,b})_{p^\infty}$ is finite, that every prime of $F$ is ramified in $F_{a,b}/F$, no prime in $S$ splits completely in $F_{a,b}/F$ and that $X(A/F_{a,b})$ is torsion over $\Zp\ps{\Gal(F_{a,b}/F)}$, then the proof of Theorem \ref{main thm Sel} carries over to yield the following inequality
      \[ \#\M(A/F_\infty)\leq \ord_{U_{a,b}}\left(\ch_{\Zp\ps{\Gal(F_{a,b}/F)}} X(A/F_{a,b})\right), \]
      where we identified $\Zp\ps{\Gal(F_{a,b}/F)}\cong \Zp\ps{U_{a,b}}$. The reason of imposing the ramification/decomposition condition is necessitated by the argument in analysing the kernel of the local maps $g_v$. As seen in the proof of Lemma \ref{control theorem}, when dealing with primes above $p$, we need to make use of results of Coates-Greenberg \cite{CG} to give an alternative description of $J_v(A/-)$, and this alternative description is necessary for us to apply a restriction-inflation sequence to analysis the kernel. In order for their result to be applicable, the ramification assumption is therefore required. On the other hand, when discussing the kernel for primes outside $p$, we need to be able to express the $J_v(A/F_{a,b})$ as a finite sum, which therefore forces us to work with the assumption that no prime in $S$ splits completely in $F_{a,b}/F$.
  \item[$(3)$]  If $F_{a,b}$ is some $\Zp$-extension of $F$ with the properties that every prime of $F$ is ramified in $F_{a,b}/F$ and no prime in $S$ splits completely in $F_{a,b}/F$, the authors of \cite{KMS} showed an upper bound in term of the $\la$-invariant of $X(A/F_{a,b})$. The bound presented here, and also that in \cite{GHKL}, is therefore an improvement of theirs.
\end{enumerate}
\er

We end the section with a supplementary result. As mentioned above, the authors of \cite{GHKL} has established Theorem \ref{main thm Sel} under the extra assumption that $X(A/F_\infty)$ is a direct sum of cyclic torsion modules. The main reason of this extra assumption is because for their proof approach, they need to show that \[ \pi\big(\ch_{\Zp\ps{G}}X(A/F_\infty)\big) = \ch_{\Zp\ps{\Ga}}X(A/F_\cyc), \]
where $\pi:\Zp\ps{G}\twoheadrightarrow\Zp\ps{\Ga}$ is the map induced by the natural projection $G\twoheadrightarrow \Ga$.
Although we do not require this for our eventual proof, we thought that it would be of interest to record the observation that this said identity can be established without the said additional assumption that $X(A/F_\infty)$ is a direct sum of cyclic torsion $\Zp\ps{G}$-modules. This is the content of the supplementary result.

\bp \label{Sel lemma}
Retain the assumptions of Theorem \ref{main thm Sel}. Then the following statements are valid.
\begin{enumerate}
    \item[$(a)$] For $K_\infty \in \{F_\cyc, F_\infty\}$, the sequence
  \[ 0\lra\Sel(A/K_\infty) \lra H^1(G_S(K_\infty),A_{p^\infty})\lra \bigoplus_{v\in S}J_v(A/K_\infty)\lra 0 \] is short exact and $H^2(G_S(K_\infty),A_{p^\infty})=0$.
  \item[$(b)$] $H_1(H,X(A/F_\infty))=0$.
  \item[$(c)$]  $\pi\big(\ch_{\Zp\ps{G}}X(A/F_\infty)\big) = \ch_{\Zp\ps{\Ga}}X(A/F_\cyc)$.
\end{enumerate}
\ep

\bpf
As seen in the proof of the main theorem, we have that $X(A/F_\infty)$ is torsion over $\Zp\ps{G}$. Next, we recall that by  the Poitou-Tate exact sequence, we have an exact sequence
\begin{multline*}
  0\lra\Sel(A/K_\infty) \lra H^1(G_S(K_\infty),A_{p^\infty})\lra \bigoplus_{v\in S}J_v(A/K_\infty)  \\
 \lra \mathfrak{S}(A^*/K_\infty)^\vee \lra H^2(G_S(K_\infty),A_{p^\infty})\lra 0,
\end{multline*}
where $\mathfrak{S}(A^*/K_\infty)$ is a submodule of $H^1_\Iw(K_\infty/F, T^*):= \plim_L H^1(G_S(L),T^*)$ with $T^*=T_pA^*$ being the Tate module of the dual abelian variety $A^*$ of $A$. Thus, the verification of (a) is reduced to proving $\mathfrak{S}(A^*/K_\infty)=0$.
Now consider the following spectral sequence of Jannsen-Nekov\'a\v{r} (\cite{jannsen, Nek})
\[ \Ext^i_{\Zp\ps{\Gal(K_\infty/F)}}\big(H^j(G_S(K_\infty),A^*_{p^\infty})^\vee,\Zp\ps{\Gal(K_\infty/F)}\big) \Rightarrow
H^{i+j}_\Iw(K_\infty/F, T^*)\]
The low degree terms fit into the following exact sequence
\[ 0\lra \Ext^1_{\Zp\ps{\Gal(K_\infty/F)}}\big(A^*_{p^\infty}(K_\infty)^\vee,\Zp\ps{\Gal(K_\infty/F)}\big)
\lra H^{1}_\Iw(K_\infty/F, T^*)\] \[\lra \Ext^0_{\Zp\ps{\Gal(K_\infty/F)}}\big(H^1(G_S(K_\infty),A^*_{p^\infty})^\vee,\Zp\ps{\Gal(K_\infty/F)}\big). \]
Now, if $K_\infty=F_\infty$, then $A^*_{p^\infty}(F_\infty)^\vee$ is plainly pseudo-null over $\Zp\ps{\Gal(K_\infty/F)}$ since in this context $\Gal(K_\infty/F)\cong \Zp^2$ and $A^*(F_\infty)_{p^\infty}^\vee$ is at most finitely generated over $\Op$.
If $K_\infty = F_\cyc$, then a result of Imai \cite{Im} tells us that $A^*_{p^\infty}(F_\cyc)^\vee$ is finite, and so, is pseudo-null over $\Zp\ps{\Gal(K_\infty/F)}$. Either way, the $\Ext^1$-term in the above exact sequence vanishes, and as a consequence, we obtain injections
\[ \mathfrak{S}(A^*/K_\infty) \hookrightarrow H^{1}_\Iw(K_\infty/F, T^*)\hookrightarrow \Ext^0_{\Zp\ps{\Gal(K_\infty/F)}}\big(H^1(G_S(K_\infty),A^*_{p^\infty})^\vee,\Zp\ps{\Gal(K_\infty/F)}\big). \]
Since an $\Ext^0$-term does not contain non-trivial torsion submodule, and so does $\mathfrak{S}(A^*/K_\infty)$.

On the other hand, taking the torsionness of $X(A/F_{\infty})$, the formulas in \cite[Theorem 4.1]{OV} into account, and followed by a straightforward rank calculation, we have that $\mathfrak{S}(A^*/K_{\infty})^{\vee}$ has zero $\Zp\ps{G}$-rank. Hence this forces $\mathfrak{S}(A^*/K_{\infty})=0$, as required.

We now verify statement (b). Taking (a) into account, we have the following diagram
\[   \entrymodifiers={!! <0pt, .8ex>+} \SelectTips{eu}{}\xymatrix{
    0 \ar[r] &\Sel(A/F_c)  \ar[d]^{\al} \ar[r] &  H^1\big(G_{S}(F_\cyc),A_{p^\infty}\big)
    \ar[d]^{h} \ar[r] & \displaystyle \bigoplus_{v\in S}J_v(A/F_\cyc) \ar[d]^{g=\oplus g_v} \ar[r] & 0  \\
    0 \ar[r]^{} & \Sel(A/F_\infty)^{H} \ar[r]^{} & H^1\big(G_{S}(F_\infty),A_{p^\infty}\big)^{H} \ar[r] &
    \left(\displaystyle \bigoplus_{v\in S}J_v(A/F_\infty)\right)^{H} \ar@{->}`r/8pt[dl]`[dl]`[ddlll]`[ddr] &\\
      &&&& \\
     & H^1(H,\Sel(A/F_\infty)) \ar[r] & H^1\big(H,H^1\big(G_{S}(F_\infty),A_{p^\infty}\big)\big)& &} \]

with exact rows. Since $\mathrm{cd}_p(H)=1$, the map $g$ is surjective. Therefore, we are reduced to showing that
\[ H^1\big(H,H^1\big(G_{S}(F_\infty),A_{p^\infty}\big)\big)=0.\]
Indeed, in view that $H^2(G_S(F_\infty),A_{p^\infty})=0$ and that $\mathrm{cd}_p(H)=1$, the Hochschild-Serre spectral sequence
\[ H^i(H, H^j(G_S(F_\infty),A_{p^\infty}))\Longrightarrow H^{i+j}(G_S(F_\cyc),A_{p^\infty}) \]
degenerates yielding
\[ H^1\big(H,H^1\big(G_{S}(F_\infty),A_{p^\infty}\big)\big)\cong H^2(G_S(F_\cyc),A_{p^\infty}),\]
 where the latter vanishes by (a). This completes the proof of (b). Assertion (c) now follows as a consequence of this and Lemma \ref{charproject}.
\epf

\section{Mazur growth number conjecture} \label{Mazur Growth section}

We now apply Theorem \ref{main thm Sel} to study the Mazur growth conjecture (see \cite[Section 18]{Maz84} for the original statement of the conjecture; also see \cite{GHKL, KMS, KL, KL2} for recent related studies). In particular, we will improve some of the results in \cite{GHKL, KL}.  To begin with, we have the following.

\bp
Let $E$ be an elliptic curve defined over $\Q$ with conductor $N$ such that it has good ordinary reduction at $p$. Let $K$ be an imaginary quadratic field of $\Q$ which satisfies all of the following.
\begin{enumerate}
  \item[$(i)$] The prime $p$ splits in $K/\Q$, and every prime divisor of $N$ is unramified in $K/\Q$.
  \item[$(ii)$] One has $N = N^+N^-$ such that $N^+$ is the largest
factor of $N$ divisible only by primes that are split in $K$, and $N^-$ is a squarefree product
of an even number of primes all of which are inert in $K$.
  \item[$(iii)$] The two primes of $K$ above $p$ are totally ramified in $K_{\mathrm{ac}}$, where $K_{\mathrm{ac}}$ is the anti-cyclotomic $\Zp$-extension of $K$.
  \item[$(iv)$] $\ord_U \Big(\ch_{\Zp\ps{\Ga}} X(E/K_\cyc) \Big)=1$.
\end{enumerate}
Then the Mordell-Weil rank of $E$  is bounded in every $\Zp$-extension of $K$, except in the case of $K_{\mathrm{ac}}$.
\ep

\bpf
From the discussion in \cite[Appendix A.1]{LV}, we see that $X(E/K_{\mathrm{ac}})$ is not torsion over $\Zp\ps{\Gal(K_{\mathrm{ac}}/K)}$ under hypotheses (i)-(iii). Hence, upon taking
Theorem \ref{main thm Sel} and hypothesis (iv) into consideration, we have the conclusion of the proposition.
\epf

\br In particular, the preceding result improves that in \cite[Corollary 3.10]{GHKL}, where we remove the restrictive assumption that the Pontryagin dual of the Selmer group of $E$ over the $\Zp^2$-extension of $K$ is a direct sum of cyclic torsion modules.
\er

We now consider an analog of the above result, which is also proven in \cite[Theorem 1.3]{KL}.
Let $\mathcal{A}$ be a simple modular self-dual abelian variety of $\mathrm{GL}_2$-type and level $\N$ over a totally real field $F$ with trivial central character. We always assume that $\mathcal{A}$ has potential good ordinary reduction at all primes above $p$. Let $K$ be a totally
imaginary extension of $F$ such that $\epsilon_{K/F}(\N)=(-1)^{|F:\Q|-1}$, where $\epsilon_{K/F}$ is the quadratic character attached to $K/F$. We further assume that every prime of $F$ above $p$ splits in $K/F$. Denote by $K_{\ac}$ the anti-cyclotomic extension of $K$. Note that $K_{\ac}$ is not necessarily a $\Zp$-extension but possibly a multiple $\Zp$-extension of $K$.

Let $\mathfrak{p}$ be a prime of $F$ above $p$. Denote by $\Op$ the ring of integer of the completion of $F$ at $\mathfrak{p}$. We then denote by $\mathrm{Sel}_{\mathfrak{p}^\infty}(\mathcal{A}/-)$ the $\mathfrak{p}$-primary Selmer group which is defined in a manner analogous to the usual $p$-primary Selmer group, except that ``$A_{p^\infty}$'' is replaced by ``$\mathcal{A}_{\mathfrak{p}^\infty}$''. The Pontryagin dual of this Selmer group will be denoted $X_{\mathfrak{p}^\infty}(\mathcal{A}/-)$. Our result in this context is as follows.

\bp
Retain the above settings. Suppose further that $X_{\mathfrak{p}^\infty}(\mathcal{A}/K_\cyc)$ is torsion over $\Op\ps{\Ga}$ with
$\ord_{U} \big(\ch_{\Op\ps{U}}X_{\mathfrak{p}^\infty}(\mathcal{A}/K_\cyc)\big)=1$. Then for every $\Zp$-extension $\mathcal{K}_\infty$ contained in $K_{\mathrm{ac}}$, the following statements are valid.
\begin{enumerate}
  \item[$(a)$] $X_{\mathfrak{p}^\infty}(\mathcal{A}/\mathcal{K}_\infty)$ has $\Op\ps{\Gal(\mathcal{K}_\infty/K)}$-rank $1$.
  \item[$(b)$] $X_{\mathfrak{p}^\infty}(\mathcal{A}/L_\infty)$ is torsion over $\Op\ps{\Gal(L_\infty/K)}$ for every $\Zp$-extension $L_\infty\neq \mathcal{K}_\infty$ contained in the compositum $K_\cyc\mathcal{K}_\infty$.
\end{enumerate}
\ep

\bpf
 By an $\mathfrak{p}$-analog of Mazur Control Theorem (see \cite{MO}), the restriction map
\[ \mathrm{Sel}_{\mathfrak{p}^\infty}(\mathcal{A}/K) \lra \mathrm{Sel}_{\mathfrak{p}^\infty}(\mathcal{A}/M_\infty)^{\Gal(M_\infty/K)}\]
has finite kernel and cokernel for every $\Zp$-extension $M_\infty$ of $K$.
Applying this to $K_\cyc$, we see that the module $\mathrm{Sel}_{\mathfrak{p}^\infty}(\mathcal{A}/K)$ has $\Op$-corank $1$. Taking this into account, and turning to $\mathcal{K}_\infty$, we have that $\mathrm{Sel}_{\mathfrak{p}^\infty}(\mathcal{A}/\mathcal{K}_\infty)^{\Gal(\mathcal{K}_\infty/K)}$ has $\Op$-corank $1$. By \cite[Proposition 5.3.20]{NSW}, one has
\begin{eqnarray*}
  \rank_{\Op\ps{\Gal(\mathcal{K}_\infty/K)}} X_{\mathfrak{p}^\infty}(\mathcal{A}/\mathcal{K}_\infty) \! &=&\! \rank_{\Op} X_{\mathfrak{p}^\infty}(\mathcal{A}/\mathcal{K}_\infty)_{\Gal(\mathcal{K}_\infty/K)}- \rank_{\Op} X_{\mathfrak{p}^\infty}(\mathcal{A}/\mathcal{K}_\infty)^{\Gal(\mathcal{K}_\infty/K)} \\
   \!&=&\! 1- \rank_{\Op} X_{\mathfrak{p}^\infty}(\mathcal{A}/\mathcal{K}_\infty)^{\Gal(\mathcal{K}_\infty/K)}.
\end{eqnarray*}
On the other hand, it follows from \cite[Theorem 0.4]{Nek08} that the module $X_{\mathfrak{p}^\infty}(\mathcal{A}/\mathcal{K}_\infty)$ is not torsion over $\Op\ps{\Gal(\mathcal{K}_\infty/K)}$, and so its $\Op\ps{\Gal(\mathcal{K}_\infty/K)}$-rank is $\geq 1$. Taking this observation and the above equality into account, we see that $X_{\mathfrak{p}^\infty}(\mathcal{A}/\mathcal{K}_\infty)$ is forced to have $\Op\ps{\Gal(\mathcal{K}_\infty/K)}$-rank $1$, and this proves assertion (a). Assertion (b) follows directly from Theorem \ref{main thm Sel}.
\epf

\br The preceding result improves that in \cite[Theorem 6.2]{KL}, where we remove the assumption that the Pontryagin dual of the Selmer group of $E$ over the $\Zp^2$-extension of $K$ is a direct sum of cyclic torsion modules.
\er

\section{Diophantine stability and integrally Diophantine extension} \label{DioStable section}

We now apply Theorem \ref{main thm Sel} to study the question of Diophantine stability and integrally Diophantine extensions. To begin with, we let $A$ be an abelian variety defined over a number field over $F$, and $L$ a finite extension of $F$. Following Mazur-Rubin \cite{MRL}, we say that the abelian variety $A$ is diophantine-stable for $L/F$ if $A(L)=A(F)$.

Let $n$ be a given positive integer. We will be concerned with the following set $\mathfrak{D}_n(A/F)$ which consists of pair $(L_1, L_2)$ satisfying the following:

 \begin{enumerate}
   \item[$(1)$] $L_1$ and $L_2$ are contained in a $\Zp$-extension of $F$ with $|L_1:L_2|=p^n$.
   \item[$(2)$] The abelian variety $A$ is diophantine-stable for $L_1/L_2$.
 \end{enumerate}

%We also define the set $\mathfrak{R}_n(A/F)$ which consists of pair $(L_1, L_2)$ satisfying the following:

% \begin{enumerate}
%   \item[$(1)$] $L_1$ and $L_2$ are contained in a $\Zp$-extension of $F$ with $|L_1:L_2|=p^n$.
%   \item[$(2)$] $\rank_\Z A(L_1) = \rank_\Z A(L_2)$.
% \end{enumerate}

We are now in position to state our result pertaining to Diophantine stability.

\bt \label{DioStable}
Let $F$ be a number field with at least one complex prime.
Suppose that $A$ is an abelian variety defined over $F$ with potential good ordinary reduction at every prime above $p$. Suppose that $X(A/F_\cyc)$ is torsion over $\Zp\ps{\Gal(F_\cyc/F)}$. For every integer $n\geq 1$, the set
$\mathfrak{D}_n(A/F)$ is uncountable.
\et

\bpf
Let $F_\infty$ be a $\Zp^2$-extension of $F$ which contains $F_\cyc$. For each $(a,b)\in \mathbb{P}_1(\Zp)\setminus\M(\mathcal{A}/F_\infty)$, the Mordell-Weil rank of $A$ is bounded in $F_{a,b}$. On the other hand, a result of Wingberg \cite[Theorem 4.3]{Win} tells us that there is at most finitely many $\Zp$-extension $K_\infty$ of $F$ such that $A(K_\infty)$ has an infinite torsion group. (Remark: Wingberg result's is stated for a simple abelian variety, but in general, an abelian variety $A$ is isogenous to a finite product of simple abelian varieties, and so this finite number of simple abelian varieties will still give finitely many many $\Zp$-extension $K_\infty$ of $F$ such that $A(K_\infty)$ has an infinite torsion group.) Hence we have uncountably many $(a,b)\in \mathbb{P}_1(\Zp)$ such that the Mordell-Weil rank of $A$ is bounded in $F_{a,b}$ and the torsion subgroup of $A(F_{a,b})$ is finite. It then follows from a classical argument (for instances, see \cite[Theorem 1.3]{GPark}) that $A(F_{a,b})$ is finitely generated for each of such $(a,b)$. Therefore, we can always find subextensions $L_1/L_2$ in $F_{a,b}$ which lies in $\mathfrak{D}_n(A/F)$.
\epf

Combining the above with the results of Rubin and Kato, we have the following unconditional result.

\bc \label{DioStable coro}
Let $E$ be an elliptic curve defined over $\Q$ with good ordinary reduction at $p$.
Suppose that $F$ is a finite abelian imaginary extension of $\Q$. Then for every integer $n\geq 1$, the set
$\mathfrak{D}_n(E/F)$ is uncountable.
\ec

We now come to the topic of integrally diophantine extensions, where we begin by recalling their definitions.

\bd
Let $R$ be a commutative ring with identity and let $n$ be a given positive
integer. We say that a subset $\Sigma$ of $R$ is Diophantine in $R$ if there exists a positive integer $m$ and a polynomial $P(x, y_1,..., y_m)$ with coefficients in $R$ such that
$a$ is in $\Sigma$ if and only if there exist elements $b_1, . . . , b_m$ of $R$ for which
$P(a, b_1, . . . , b_m) = 0$.

Let $L_1/L_2$ be an extension of number fields. If the ring of integers of $L_2$ is Diophantine in the ring of integers of $L_1$, then $L_1/L_2$ is said to be integrally Diophantine.
\ed

A conjecture of Denef and Lipshitz predicts that $K/\Q$ is integrally diophantine for every number field $K$ (see \cite{DL}), and the validity of this conjecture has been established by Koyman and Pagano (see \cite{KP}). It's then natural ask whether the following will always hold.

\medskip \noindent \textbf{Relative Denef-Lipshitz conjecture.}
Every finite extension of number field $L/K$ is integrally diophantine.

\medskip
The following theorem of Shlapentokh \cite[Theorem 1.9]{Sh} (also see \cite[Theorem 3.1]{MRS}) is a fundamental tool for the study of these type of problems.

\bt[Shlapentokh]
Let $L_2/L_1$ be an extension of number fields. Suppose that
there is an abelian abelian $A$ defined over $L_1$ such that $\rank_\Z A(L_2) = \rank_\Z A(L_1) > 0$. Then $L_2/L_1$ is integrally Diophantine.
\et

Let $\mathcal{I}_n(F)$ be the set which consists of pair $(L_1, L_2)$ satisfying the following:

 \begin{enumerate}
   \item[$(1)$] $L_1$ and $L_2$ are contained in a $\Zp$-extension of $F$ with $|L_1:L_2|=p^n$.
   \item[$(2)$] $L_1/L_2$ is integrally Diophantine.
 \end{enumerate}

 Our result concerning integrally Diophantine extensions is as follows.

\bt \label{intDio}
Let $F$ be a number field with at least one complex prime.
Suppose that there exist an abelian variety $A$ defined over $F$ which satisfies all of the following.
\begin{enumerate}
  \item[$(i)$] The abelian variety $A$ has potential good ordinary reduction at every prime above $p$.
  \item[$(ii)$] $X(A/F_\cyc)$ is torsion over $\Zp\ps{\Gal(F_\cyc/F)}$.
  \item[$(iii)$] $\rank_\Z A(F) >0$.
\end{enumerate}
Then, for every integer $n\geq 1$, the set
$\mathcal{I}_n(F)$ is uncountable.
\et

\bpf
As seen in the proof of Theorem \ref{DioStable}, there exists uncountably many $\Zp$-extensions of $F$ such that the Mordell-Weil rank of $A$ is bounded. Therefore, for each such $\Zp$-extension, one can always find intermediate extension $L_1/L_2$ with $\rank_\Z A(L_1) = \rank A(L_2)$, and this common rank is necessarily greater than $0$ by hypothesis (iii). Hence we may apply Shlapentokh's theorem to conclude that $L_1/L_2$ is integrally Diophantine.
\epf

Unconditionally, we have the following.

\bc \label{intDio coro}
Suppose that $F$ is a finite abelian imaginary extension of $\Q$. Then for every integer $n\geq 1$, the set
$\mathcal{I}_n(F)$ is uncountable.
\ec

\bpf
By \cite[Proposition 5.4 and Remark after which]{G99}, there exist an elliptic curve $E$ defined over $\Q$ having good ordinary reduction at $p$ with $L(E/\Q,1)\neq 0$. By appealing to the work of Bump-Friedberg-Hoffstein \cite{BFH} or that of Murty-Murty \cite{MM}, we can find a quadratic twist $E^{(D)}$ of $E$ such that $L(E^{(D)}/\Q,s)$ has a simple zero at $s=1$. From the deep result of Gross-Zagier \cite{GZ}, this in turn implies $\rank_{\Z}E^{(D)}(\Q)=1$. In particular, we have $\rank_{\Z}E^{(D)}(F)>0$. On the other hand, a result of Kato \cite{K} tells us that $X(E^{(D)}/F_\cyc)$ is torsion over $\Zp\ps{\Gal(F_\cyc/F)}$. Therefore, all the hypothesis of Theorem \ref{intDio} are satisfied, and this gives the conclusion of the corollary.
\epf

We end the section with another result which is useful in obtaining cases of uncountable $\mathfrak{D}_n(L)$ and $\mathcal{I}_n(L)$ for non-abelian $L$.

\bp \label{Kida}
Let $F$ be a number field with at least one complex prime.
Suppose that there exist an elliptic curve $E$ defined over $F$ which satisfies all of the following.
\begin{enumerate}
  \item[$(i)$] The elliptic curve $E$ has potential good ordinary reduction at every prime above $p$.
  \item[$(ii)$] $X(E/F_\cyc)$ is finitely generated over $\Zp$.
\end{enumerate}
Then for every finite Galois $p$-extension $L$ of $F$ and every integer $n\geq 1$, the set
$\mathfrak{D}_n(E/L)$ is uncountable.

In the event that  $\rank_\Z E(F) >0$, we even have that the set
$\mathcal{I}_n(L)$ is uncountable.
\ep

\bpf
It's well-known that under assumptions of the proposition, $X(E/L_\cyc)$ is finitely generated over $\Zp$ for every finite Galois $p$-extension $L$ of $F$ (for instance, see \cite[Corollary 3.4]{HM}). Therefore, in particularly, $X(E/L_\cyc)$ is a torsion $\Zp\ps{\Gal(L_\cyc/L)}$-extension. Therefore, we may apply Theorems \ref{DioStable} and  \ref{intDio} to obtain the conclusion of the proposition.
\epf

We give an example to illustrate Proposition \ref{Kida}. Let $p=3$ and let $E$ be the elliptic curve 79A1 of Cremona tables. It follows from the discussion in \cite[Page 253]{DD} that $\rank_\Z E(\Q(\mu_3))=1$ and $X(E/\Q(\mu_{3^\infty}))\cong \Zp$. Therefore, Proposition \ref{Kida} tells us that $\mathfrak{D}_n(E/L)$ and $\mathcal{I}_n(L)$ are uncountable for every finite Galois $3$-extension $L$ of $\Q(\mu_3)$. In particular, for instance, the set  $\mathcal{I}_n(\Q(\mu_{3^n}, \sqrt[3^n]{m}))$ is infinite for every cubefree integer $m$ and positive integer $n$.

\section{Some further remark on Theorem \ref{main thm Sel}} \label{Remark on Thm3.3 section}

We retain the notation and setting of Section \ref{main thm section}. In particular, we continue to identify $\Zp\ps{\Ga}$ with $\Zp\ps{U}$, where $\Ga= \Gal(F_\cyc/F)$, and we continue to assume that $X(A/F^\cyc)$ is torsion over $\Zp\ps{\Ga}$. Then the structure theory of $\Zp\ps{\Ga}$-module tells us that there is a pseudo-isomorphism
\[X(A/F^\cyc) \sim \bigoplus_{i=1}^s\Zp\ps{\Ga}/p^{\al_i} \times \bigoplus_{j=1}^t\Zp\ps{\Ga}/f_j^{\be_j},\]
where each $f_j$ is irreducible in $\Zp\ps{\Ga}$ and is not an associate of $p$.

If one assumes the elementary factors of $X(A/F^\cyc)$ do not contain $\Zp\ps{U}/U^e$ for every $e\geq 2$, it then follows from Mazur's control theorem that
 \[\ord_{U}\Big(\ch\big(X(\mathcal{A}/F_\cyc)\big)\Big) = \corank_{\Zp}\big(\Sel(\mathcal{A}/F)\big).\]
 We note that the above assumption is plainly true if $\ord_{U}\Big(\ch\big(X(\mathcal{A}/F_\cyc)\big)\Big)\leq 1$.
 (for instance, see \cite[Page 59, 1st paragraph]{G99} or \cite[Proposition 4.5]{LimVan}).
In fact, this assumption is a consequence of the following semi-simplicity conjecture of Greenberg \cite[Conjecture 1.12]{G99}.

\begin{conjecture}[Greenberg] \label{semisimple conj}
$\be_j =1$ for every $j$.
\end{conjecture}

Plainly, this conjecture of Greenberg is true if the Iwasawa $\lambda$-invariant is $1$. To the best knowledge of the author, there seems very little evidence in literature on Conjecture \ref{semisimple conj}. In \cite[Lemma 4.6]{LimVan}, the author gives a simple criterion (which is far from being sufficient!) for verifying this conjecture.

Granted this conjecture of Greenberg, we have the following interesting (conjectural) observation.

\bp \label{general ab var}
Let $A$ be an abelian variety defined over $F$ with potential good ordinary reduction at every prime above $p$. Let $F_\infty$ be a $\Zp^2$-extension of $F$ which contained the cyclotomic $\Zp$-extension $F_\cyc$. Suppose that $X(A/F_\cyc)$ is torsion over $\Zp\ps{\Gal(F_\cyc/F)}$. If we assume that the semisimplicity conjecture of Greenberg holds, then we have
\[ \#\M(A/F_\infty)\leq  \corank_{\Zp}\big(\Sel(A/F)\big). \]
If we further assume that the $p$-primary Tate-Shafarevich group of $A$ over $F$ is finite, then one has the inequality
\[ \#\M(A/F_\infty)\leq \rank_{\Z} A(F). \]
\ep

The above observation, though conjectural, is very interesting, as it is saying that the number of $\Zp$-extensions in a $\Zp^2$-extension at which $A$ has unbounded growth is bounded above by the Mordell-Weil rank of $A$ at the base field! It is also interesting to note that this upper bound is independent of the $\Zp^2$-extension (as long as it contains $F_\cyc$) in question.

\section{Signed Selmer groups of elliptic curve with supersingular reduction} \label{elliptic supersingular section}

We come to the final section of the paper which will study an elliptic curve with good supersingular reduction at $p$. We first consider the local situation. Suppose for now that $E$ is an elliptic curve defined over $\Qp$ which has good supersingular reduction.  We shall also assume that $a_p= p+1 - |\widetilde{E}(\mathbb{F}_p)|=0$, where $\widetilde{E}$ is the reduced curve of $E$ modulo $p$. We write $\widehat{E}$ for the formal group of $E$.

Let $k_m$ denote the unique unramified extension $\Qp$ such that $|k_m:\Qp|=p^m$. By convention, we set $k_0=\Qp$. Then $k_\infty= \cup_{n\geq 0}k_n$ is the unramified $\Zp$-extension of $\Qp$. Write $\Q_{p,n}$ the unique subextension of $\Q_{p,\cyc}/\Qp$ with $|\Q_{p,n}:\Qp|=p^n$. Let $L_\infty$ be another $\Zp$-extension of $\Qp$ which is neither $\Q_{p,\cyc}$ nor $k_\infty$. Such a $\Zp$-extension is necessarily totally ramified, and we let $L_n$ be the subextension of $L_\infty/\Qp$ with $|L_n:\Qp|=p^n$. As before, we take $L_0=\Qp$.

\bl \label{local fields agree}
For $m\geq n\geq 0$, we have $k_mL_n = k_m\Q_{p,n}$.
\el

\bpf (cf. \cite[Lemma 2.13]{Kim})
Since $L_\infty(\mu_p)$ is a $\Zp^\times$-extension of $\Qp$, its group of universal norms is generated by a uniformizer $\varpi$ of $\Zp$ such that $\ord_p(\frac{\varpi}{p}-1)>0$. By local class field theory, we have
\[\Gal(\Q_{p,n}/\Qp)\cong \Qp^{\times}/\langle p \rangle\cdot \mu(\Qp^\times)\cdot(1+p^{n+1}\Zp) \]
and
\[ \Gal(L_n/\Qp)\cong \Qp^{\times}/\langle \varpi\rangle\cdot \mu(\Qp^\times)\cdot(1+p^{n+1}\Zp)  \]
Since $k_mL_n$  (resp. $k_m\Q_{p,n}$) is an unramified extension of $L_n$ (resp., of $\Q_{p,n}$) with degree $p^m$, we have
\begin{eqnarray*}
  \Gal(k_m\Q_{p,n}/\Qp)\cong \Qp^{\times}/\langle p^{p^m} \rangle\cdot \mu(\Qp^\times)\cdot(1+p^{n+1}\Zp) \\
 \Big(\mbox{resp.,}  \Gal(k_mL_n/\Qp)\cong \Qp^{\times}/\langle \varpi^{p^m} \rangle\cdot \mu(\Qp^\times)\cdot(1+p^{n+1}\Zp) \Big).
\end{eqnarray*}
Now since we have  $\ord_p(\frac{\varpi}{p}-1)>0$, it follows that
\[ \left(\frac{\varpi}{p}\right)^{p^m}\in 1+p^{m+1}\Zp \subseteq 1+p^{n+1}\Zp,\]
where the second inclusion is a consequence of our assumption that $m\geq n$. The conclusion of the lemma then follows from these observations.
\epf

\bd \label{local norm groups} Following \cite{IP, Kim, Kob}, we define the following plus/minus norm groups.
\[\widehat{E}^+(k_mL_n) = \{ P\in \widehat{E}(k_mL_n)~:~\mathrm{tr}_{n/l+1,m}(P)\in \widehat{E}(k_mL_l), 2\mid l, l<n \}, \]
\[\widehat{E}^-(k_mL_n) = \{ P\in \widehat{E}(k_mL_n)~:~\mathrm{tr}_{n/l+1,m}(P)\in \widehat{E}(k_mL_l), 2\nmid l, l<n \}, \]
where $\mathrm{tr}_{n/l+1,m}:\widehat{E}(k_mL_n) \lra \widehat{E}(k_mL_{l+1})$ denotes the trace map. In the event that $m=0$, we shall write the plus/minus norm groups as $E^\pm(L_n)$.
In the event that $n=0$, we have $\widehat{E}^\pm(L_nk_m) = \widehat{E}^\pm(k_m) = \widehat{E}(k_m)$. It is a straightforward exercise to check that $\widehat{E}^\pm(L_n)\subseteq \widehat{E}^\pm(k_mL_n)$ and $\widehat{E}(k_m)=\widehat{E}^\pm(k_m)\subseteq \widehat{E}^\pm(k_mL_n)$

Finally, the plus/minus norm groups $\widehat{E}^+(k_m\Q_{p,n})$ and $\widehat{E}^-(k_m\Q_{p,n})$ are defined similarly.
\ed

We then define $\widehat{E}^\pm(L_\infty) = \cup_{n\geq 0} \widehat{E}^\pm(L_n)$ and $\widehat{E}^\pm(\Q_{p,\cyc}) = \cup_{n\geq 0} \widehat{E}^\pm(\Q_{p,n})$. We also write $\widehat{E}^\pm(k_\infty) = \cup_{n\geq 0}\widehat{E}^\pm(k_n)= \cup_{n\geq 0}\widehat{E}(k_m)$. To summarize, we have given a definition of $\widehat{E}^\pm(M_\infty)$ for every $\Zp$-extension $M_\infty$ of $F$.

Now, let $K_\infty$ be the $\Zp^2$-extension of $\Qp$. We like to give a definition for $\widehat{E}^\pm(K_\infty)$. The natural one is to take $\cup_{m,n}\widehat{E}^\pm(k_m\Q_{p,n})$ as a definition. However, for our purposes, we need to also consider taking $\cup_{m,n}\widehat{E}^\pm(k_mL_n)$ for a ramified $\Zp$-extension $L_\infty$ of $\Qp$. Thankfully, the next result tells us that the choice of the ramified $\Zp$-extension $L_\infty$ does not matter.

\bl \label{plus/minus norm group independent}
For $m\geq n\geq 0$ and $s\in\{+,-\}$, we have $\widehat{E}^s(k_mL_n) = \widehat{E}^s(k_m\Q_{p,n})$. In particular, for $m\geq n\geq 0$, the definition of the plus/minus norm group $\widehat{E}^s(k_mL_n)$ is independent of the choice of the ramified $\Zp$-extension $L_\infty$.
\el

\bpf
This is immediate from Lemma \ref{local fields agree} and the definition of the plus/minus norm groups.
\epf

\bl \label{norm groups inject}
For every $\Zp$-extension $M_\infty$ of $\Qp$, we have a map
\[ \frac{H^1(M_\infty, E_{p^\infty})}{\widehat{E}^\pm(M_\infty)\ot\Qp/\Zp} \lra \left(\frac{H^1(K_\infty, E_{p^\infty})}{\widehat{E}^\pm(K_\infty)\ot\Qp/\Zp}\right)^{\Gal(K_\infty/M_\infty)}. \]
In the event that $M_\infty\neq k_\infty$, the above map is an isomorphism.
\el

\bpf
Foe every $n$, we have the following commutative diagram
 \[   \entrymodifiers={!! <0pt, .8ex>+} \SelectTips{eu}{}\xymatrix @C=0.8pc{
    0 \ar[r] &\widehat{E}^\pm(L_n)\ot\Qp/\Zp  \ar[d]^{\al} \ar[r] &  H^1(L_n, E_{p^\infty})
    \ar[d]^{\beta} \ar[r] & \frac{H^1(L_n, E_{p^\infty})}{\widehat{E}^\pm(L_n)\ot\Qp/\Zp} \ar[d]^{\gamma}  \ar[r] &0\\
    0 \ar[r]^{} & (\widehat{E}^\pm(k_nL_n)\ot\Qp/\Zp)^{\Gal(k_nL_n/L_n)} \ar[r]^{} & H^1(k_nL_n, E_{p^\infty})^{\Gal(k_nL_n/L_n)} \ar[r] &
    \left(\frac{H^1(k_nL_n, E_{p^\infty})}{\widehat{E}^\pm(k_nL_n)\ot\Qp/\Zp}\right)^{\Gal(k_nL_n/L_n)}   & } \]
    where the injectivity in the leftmost of both rows follows from \cite[Lemma 8.17]{Kob}. The middle vertical map is the usual restriction in cohomology, the leftmost vertical map is induced by the inclusion map, and whence the leftmost square commutes which in turn induces the rightmost vertical map. Taking limit of these vertical maps, we obtain the required map of the lemma. If $M_\infty=\Q_{p,\cyc}$, this map is an isomorphism by \cite[Proposition 3.8]{LeiLim}. Upon reviewing the proof, one sees that the argument carries over to a totally ramified $\Zp$-extension of $\Qp$.
\epf

\bd \label{local norm groups ZP2}
  Let $L_\infty$ be a ramified $\Zp$-extension of $\Qp$. We define $\widehat{E}(K_\infty)=\cup_{m,n}\widehat{E}^\pm(k_mL_n)$. By cofinality, it follows from this definition that $\widehat{E}(K_\infty)=\cup_{n}\widehat{E}^\pm(k_nL_n)$. Taking this observation and Lemma \ref{norm groups inject} into account, we see that this definition is independent of the choice of the ramified $\Zp$-extension $L_\infty$.
\ed

We turn to the global situation.
From now on, $E$ will denote an elliptic curve over $\Q$ which has good supersingular reduction at the prime $p$. Denote by $\widetilde{E}$ the reduced curve of $E$ modulo $p$. We shall assume that $a_p= p+1 - |\widetilde{E}(\mathbb{F}_p)|=0$ (note that this automatically holds if $p\geq 5$). Let $F$ be an imaginary quadratic field of $\Q$ at which the prime $p$ splits completely, say $p=\p\overline{\p}$. We let $S$ be the set of primes of $F$ consisting precisely of those dividing $pN$ and the infinite primes.

We further assume that the primes $\p$ and $\overline{\p}$ are totally ramified in $F_{\ac}/F$, where $F_\ac$ is the anticyclotomic $\Zp$-extension of $F$. The $\Zp^2$-extension of $F$ is denoted by $F_\infty$. Write $F(\p^{\infty})$ for the unique
$\Zp$-extension of $F$ unramified outside $\p$ and write $F(\p^n)$ for the intermediate subfield of $F(\p^{\infty})$ with $|F(\p^n):F| =p^n$. We have analogous definitions for $F(\overline{\p}^{\infty})$ and $F(\overline{\p}^n)$. For each pair of nonnegative integers $m$ and $n$, write $F(\p^m\overline{\p}^n)$ for the compositum of the fields $F(\p^m)$ and $F(\overline{\p}^n)$.

Let $\mathcal{L}_\infty$ be a $\Zp$-extension of $F$ which is not equal to $F(\p^{\infty})$ or $F(\overline{\p}^{\infty})$. Therefore, every prime of $F$ above $p$ must ramify in $\mathcal{L}_\infty/F$. Let $s,z\in\{+,-\}$. If we write $\mathcal{L}_n$ for the intermediate subextension of $\mathcal{L}_\infty$, the signed Selmer group of $E$ over $\mathcal{L}_n$ is defined to be
\begin{multline*}
  \mathrm{Sel}^{s,z}(E/\mathcal{L}_n)=\ker\Bigg(H^1(G_S(\mathcal{L}_n)), E_{p^\infty}) \lra  \Bigg(\bigoplus_{w\mid \p}\frac{H^1( \mathcal{L}_{n,w},E_{p^\infty})}{\widehat{E}^s(\mathcal{L}_{n,w})\ot\Qp/\Zp} \Bigg) \\
  \oplus \Bigg(\bigoplus_{\overline{w}\mid \overline{\p}}\frac{H^1(\mathcal{L}_{n,\overline{w}},E_{p^\infty})}
  {\widehat{E}^z(\mathcal{L}_{n,\overline{w}})\ot\Qp/\Zp} \Bigg) \oplus \bigoplus_{u|N}H^1(\mathcal{L}_{n,u},E_{p^\infty})\Bigg).
\end{multline*}
Here $\widehat{E}^s(\mathcal{L}_{n,w})$ and $\widehat{E}^z(\mathcal{L}_{n,\overline{w}})$ are in sense as in Definition \ref{local norm groups}, where we note that $\mathcal{L}_{n,w}$ and $\mathcal{L}_{n,\overline{w}}$ are now finite extensions of $\Qp$ contained in either a ramified $\Zp$-extension or the unramified $\Zp$-extension. We then set  $\mathrm{Sel}^{s,z}(E/\mathcal{L}_\infty)=\ilim_n \mathrm{Sel}^{s,z}(E/\mathcal{L}_n)$ and write $X^{s,z}(E/\mathcal{L}_\infty)$ for the Pontryagin dual of $\mathrm{Sel}^{s,z}(E/\mathcal{L}_\infty)$.

Now let $F_n$ be the subextension of $F_\infty/F$ with $\Gal(F_n/F)\cong \Z/p^n\Z\times\Z/p\Z$. The signed Selmer group of $E$ over $F_n$ (cf. \cite{Kim}) is defined to be
\begin{multline*}
  \mathrm{Sel}^{s,z}(E/F_n))=\ker\Bigg(H^1(G_S(F_n), E_{p^\infty}) \lra \Bigg(\bigoplus_{w\mid \p}\frac{H^1(F_{n,w},E_{p^\infty}))}{E^s(F_{n,w})\ot\Qp/\Zp} \Bigg) \\
  \oplus \Bigg(\bigoplus_{\overline{w}\mid \overline{\p}}\frac{H^1(F_{n,\overline{w}},E_{p^\infty})}
  {E^z(F_{n,\overline{w}})\ot\Qp/\Zp} \Bigg) \oplus \bigoplus_{u|N}H^1(F_{n,u},E_{p^\infty})\Bigg),
\end{multline*}
Set $\mathrm{Sel}^{s,z}(E/F_\infty)=\ilim_n \mathrm{Sel}^{s,z}(E/F_n)$. We then write $X^{s,z}(E/F_\infty)$ for the Pontryagin dual of $\mathrm{Sel}^{s,z}(E/F_\infty)$.

We now make the following remark.

\br
In \cite{Kim}, Kim introduced the following groups
\[\widehat{E}^+(F(\p^m\overline{\p}^n)_w) = \{ P\in \widehat{E}(F(\p^m\overline{\p}^n)_w)~:~\mathrm{tr}_{m/l+1,n}(P)\in \widehat{E}(F(\p^l\overline{\p}^n)_w), 2\mid l, l<m \}, \]
\[\widehat{E}^-(F(\p^m\overline{\p}^n)_w) = \{ P\in \widehat{E}(F(\p^m\overline{\p}^n)_w)~:~\mathrm{tr}_{m/l+1,n}(P)\in \widehat{E}(F(\p^l\overline{p}^n)_w), 2\nmid l, l<m \}, \]
where $\mathrm{tr}_{m/l+1,n}:\widehat{E}(F(\p^m\overline{\p}^n)_w) \lra \widehat{E}(F(\p^{l+1}\overline{\p}^n)_w)$ denotes the trace map.
For a prime $\overline{w}$ of $F_{\infty}$ above $\overline{\p}$, the groups
$\widehat{E}^{\pm}(F(\p^m\overline{\p}^n)_{\overline{w}})$ are defined in a similar fashion as above. For $s,z\in\{+,-\}$, Kim defined his signed Selmer group of $E$ over $F_\infty$ by taking limit of the signed Selmer group over $F(\p^m\overline{\p}^n)$ which are given by
\begin{multline*}
  \mathrm{Sel}^{s,z}(E/F(\p^m\overline{\p}^n))=\ker\Bigg(H^1(G_S(F(\p^m\overline{\p}^n)), E_{p^\infty}) \lra  \\
  \Bigg(\bigoplus_{w\mid \p}\frac{H^1(F(\p^m\overline{\p}^n)_{w},E(p))}{\widehat{E}^s(F(\p^m\overline{\p}^n)_{w})\ot\Qp/\Zp} \Bigg)\oplus \Bigg(\bigoplus_{\overline{w}\mid \overline{\p}}\frac{H^1(F(\p^m\overline{\p}^n)_{\overline{w}},E(p))}
  {\widehat{E}^z(F(\p^m\overline{\p}^n)_{\overline{w}})\ot\Qp/\Zp} \Bigg) \oplus \bigoplus_{u|N}H^1(F(\p^m\overline{\p}^n)_u,E_{p^\infty})\Bigg),
\end{multline*}
where $S$ here is the set of primes of $F$ consisting precisely of those dividing $pN$ and the infinite primes. Taking \cite[Lemma 2.14]{Kim} into account, our signed Selmer groups over $F_\infty$ will agree with that of Kim.
\er

The following is a natural signed analogue of Mazur conjecture.

\medskip \noindent \textbf{Conjecture.} Let $F_\cyc$ be the cyclotomic $\Zp$-extension of $F$. Then $X^{s,z}(E/F_{\cyc})$ is a torsion $\Zp\ps{\Gal(F_\cyc/F)}$-module for every $s,z\in\{+,-\}$.

\medskip
When $s=z$, the conjecture is known to hold by \cite[Proposition 8.4]{LeiP} (also see \cite[Theorem 1.2]{Kob}). Note that if $X^{s,z}(E/F_\cyc)$ is torsion over $\Zp\ps{\Gal(F_\cyc/F)}$, then $X^{s,z}(E/F_\infty)$ is torsion over $\Zp\ps{G}$ (see \cite[Proposition 4.13]{LeiLim}), where $G=\Gal(F_\infty/F)$. In the subsequent discussion, we shall identity $\Zp\ps{G}$ with $\Zp\ps{U,W}$ in a way such that $\Zp\ps{W}\cong \Zp\ps{\Gal(F_\infty/F_\cyc)}$ and $\Zp\ps{U}\cong \Zp\ps{\Gal(F_\infty/F_\ac)}$. Under this choice of identification, we also have the identifications $\Zp\ps{W}\cong \Zp\ps{\Gal(F_\cyc/F)}$ and $\Zp\ps{U}\cong \Zp\ps{\Gal(F_\ac/F)}$. We continue to write $\Ga$ for $\Gal(F_\cyc/F)$. We sometimes also write $\Ga_{\ac}$ for $\Gal(F_\ac/F)$.

\bd
For $s,z\in\{+,-\}$, we define the set
\[\mathfrak{Z}^{s,z}(E/F_\infty) = \{(a,b)\in\mathbb{P}_1(\Zp)~|~\mbox{$X^{s,z}(E/F_{a,b})$ is not torsion over $\Zp\ps{\Gal(F_{a,b}/F)}$.}\} \]
\ed

We now present the following signed analogue of Theorem \ref{main thm Sel}.

\bt \label{signed X leq ord}
Let $E$ be an elliptic curve of conductor $N$ over $\Q$ with good supersingular reduction at the prime $p$ and $a_p=0$. Let $F$ be an imaginary quadratic field of $\Q$ at which all the prime divisors of $p$ split completely in $F/\Q$. If $X^{s,z}(E/F_\cyc)$ is torsion over $\Zp\ps{\Ga}$, then $X^{s,z}(E/F_\infty)$ is torsion over $\Zp\ps{G}$, and we have the following estimate
\[ \#\mathfrak{Z}^{s,z}(E/F_\infty) \leq
 \ord_U\big(\ch_{\Zp\ps{\Ga}}X^{s,z}(E/F_\cyc)\big).
\]
In particular, if $X^{s,z}(E/F_\cyc)$ is torsion over $\Zp\ps{\Ga}$ for every $s,z\in\{+,-\}$, we have
\[ \#\mathfrak{M}(E/F_\infty) \leq
 \sum_{s,z\in\{+,-\}}\ord_U\big(\ch_{\Zp\ps{\Ga}}X^{s,z}(E/F_\cyc)\big).
\]
\et

\bpf
The proof proceeds quite similarly to that in Theorem \ref{main thm Sel}, and we give a sketch of it here. We begin noting that in view of Lemma \ref{norm groups inject}, for every $(a,b)\in\mathbb{P}_1(\Zp)$, there is a restriction map
\[ \mathrm{Sel}^{s,z}(E/F_{a,b}) \lra \mathrm{Sel}^{s,z}(E/F_\infty)^{\Gal(F_\infty/F_{a,b})}\]
on the signed Selmer groups. Furthermore, this restriction map has kernel contained in $H^1(H_{a,b}, E_{p^\infty}(F_\infty))$. Since $E_{p^\infty}(F_\infty)=0$ by \cite[Proposition 8.7]{Kob}, the restriction map is therefore an injection. Therefore, it follows that if $X^{s,z}(E/F_{a,b})$ is not torsion over $\Zp\ps{\Gal(F_{a,b}/F)}$, then so is $X^{s,z}(E/F_{\infty})_{H_{a,b}}$. Thus we have the inclusion
\[ \mathfrak{Z}^{s,z}\big(E/F_\infty\big) \subseteq \mathfrak{A}\big(X^{s,z}(E/F_\infty)\big). \]
On the other hand, in the case of $F_\cyc$, the argument of the proof of \cite[Theorem 5.1]{LeiLim} shows that the restriction map
\[ \mathrm{Sel}^{s,z}(E/F_\cyc) \lra \mathrm{Sel}^{s,z}(E/F_\infty)^{\Gal(F_\infty/F_\cyc)}\]
 is bijective. Hence we may proceed similarly to that in Theorem \ref{main thm Sel} to obtain the estimate
 \[ \#\mathfrak{Z}^{s,z}(E/F_\infty) \leq
 \ord_U\big(\ch_{\Zp\ps{\Ga}}X^{s,z}(E/F_\cyc)\big).
\]
On the other hand, it follows from \cite[Corollary 4.6]{GHKL} that
\[ \mathfrak{M}(E/F_\infty)\subseteq \bigcup_{s,z\in \{+,-\}}\mathfrak{Z}^{s,z}(E/F_\infty). \]
Combining the above observations, we have the final estimate of the theorem.
\epf

\br\begin{enumerate}
   \item[$(1)$] We note that this result refines \cite[Proposition A.10]{GHKL}, where we remove the assumption that $X^{s,z}(E/F_\infty)$ is a sum of cyclic torsion modules. This assumption is required for the authors for proving
\[ \pi\big(\ch_{\Zp\ps{G}}X^{s,z}(A/F_\infty)\big) = \ch_{\Zp\ps{\Ga}}X^{s,z}(A/F_\cyc), \]
where $\pi:\Zp\ps{G}\twoheadrightarrow\Zp\ps{\Ga}$ is the map induced by the natural projection $G\twoheadrightarrow \Ga$. Although we do not require this, we mention quickly how this identity can actually be proven unconditionally. By the torsionness assumption of the signed Selmer groups, it follows from \cite[Propositions 4.6 and 4.12]{LeiLim} that we have the surjectivity of the defining sequence of the signed Selmer groups and the vanishing of $H^2(G_S(\mathcal{L}),E_{p^\infty})$ for $\mathcal{L}=F_\cyc, F_\infty$. Therefore, we may apply the argument of Proposition \ref{Sel lemma} to obtain the vanishing of $H^1(H,\mathrm{Sel}^{s,z}(E/F_\infty))$. The required identity then follows from this and Lemma \ref{charproject}.

 \item[$(2)$] In \cite[Theorem 2.18]{KMS2}, the authors obtain an upper bound in term of the sum of the $\la$-invariants of the signed Selmer groups. The bound presented here, and also that in \cite{GHKL}, is therefore an improvement of theirs. \end{enumerate}
 \er

We now apply Theorem \ref{supersingular thm} to the study of Diophantine stability for an elliptic curve with supersingular reduction at $p$. Indeed, it has been shown in \cite[Corollary 4.6]{GHKL} that if $(a,b)\in\mathbb{P}_1(\Zp)$ has the property that $X^{s,z}(E/F_{a,b})$ is torsion over $\Zp\ps{\Gal(F_{a,b}/F)}$ for every $s,z\in\{+,-\}$, then the Mordell-Weil ranks of $E$ are bounded in $F_{a,b}/F$. Therefore, if we have that
$X^{s,z}(E/F_\cyc)$ is torsion over $\Zp\ps{\Gal(F_\cyc/F)}$ for every $s,z\in\{+,-\}$, then   the set $\mathfrak{D}_n(E/F)$ is uncountable for every integer $n\geq 1$. Instead of writing this down formally, we prefer to present the following unconditional result.

\bt \label{supersingular thm}
Let $E$ be an elliptic curve of conductor $N$ over $\Q$ with good supersingular reduction at the prime $p$ and $a_p=0$. Let $F$ be an imaginary quadratic field of $\Q$ such that  $p$ splits completely in $F/\Q$, and the primes of $F$ above $p$ are totally ramified in $F_{\ac}/F$. Suppose either of the following statement is valid.
\begin{enumerate}
  \item[$(i)$] $\Sel(E/F)$ is finite.
  \item[$(ii)$] Every prime divisor of $N$ splits completely in $F/\Q$.
  \end{enumerate}
Then, for every integer $n\geq 1$, the set
$\mathfrak{D}_n(E/F)$ is uncountable.
\et

\bpf
Under assumption (i), it has been proven in \cite[Proposition 4.7]{GHKL} that $X^{s,z}(E/F_{a,b})$ is torsion over $\Zp\ps{\Gal(F_{a,b}/F)}$ for every $s,z\in\{+,-\}$ and every $(a,b)\in\mathbb{P}_1(\Zp)$. We will therefore left with establishing the theorem under assumption (ii). As seen before, the module $X^{s,s}(E/F_\cyc)$ is torsion over $\Zp\ps{\Gal(F_\cyc/F)}$ by \cite[Proposition 8.4]{LeiP}, and so Theorem \ref{signed X leq ord} tells us that $\mathfrak{Z}^{s,s}(E/F_\infty)$ is finite for $s\in\{+,-\}$. It remains to show that $\mathfrak{Z}^{s,z}(E/F_\infty)$ is finite for $s\neq z$. For this, we will go from the anticyclotomic direction. Under assumption (ii), it is shown in \cite[Theorem A.13]{GHKL} that $X^{s,z}(E/F_{ac})$ is torsion over $\Zp\ps{\Gal(F_{\ac}/F)}$. We may then apply the argument in \cite[Theorem 5.1]{LeiLim} to show that the restriction map
\[ \mathrm{Sel}^{s,z}(E/F_{\ac}) \lra \mathrm{Sel}^{s,z}(E/F_\infty)^{\Gal(F_\infty/F_{\ac})}\]
is an isomorphism, thanks to our ramification assumption on $p$ and that every prime divisor of $N$ splits in $F/\Q$. Indeed, these two assumptions ensure that the local result obtained in \cite[Proposition 3.8]{LeiLim} applies (also see Lemma \ref{norm groups inject}), and that every prime divisor of $N$ in $F$ are finitely decomposed in $F_{\ac}/F$ which is required for $J_v(E/F_{\ac})$ to be written as a finite sum (see Remark \ref{main thm Sel remark}). As a consequence, we can show that $\mathfrak{Z}^{s,z}(E/F_\infty)$ is finite with
\[ \#\mathfrak{Z}^{s,z}(E/F_\infty) \leq \ord_W\big(\ch_{\Zp\ps{\Gal(F_{ac}/F)}}X^{s,z}(E/F_{\ac})\big).
\]
The conclusion of the theorem then follows.
\epf

\br
In the setting of Theorem \ref{supersingular thm}(ii), our proof actually shows that
\begin{multline*}
  \#\mathfrak{M}(E/F_\infty) \leq
 \ord_U\big(\ch_{\Zp\ps{\Ga}}X^{+,+}(E/F_\cyc)\big) + \ord_U\big(\ch_{\Zp\ps{\Ga}}X^{-,-}(E/F_\cyc)\big) \\
 + \ord_W\big(\ch_{\Zp\ps{\Ga_{\ac}}}X^{+,-}(E/F_\ac)\big)+ \ord_W\big(\ch_{\Zp\ps{\Ga_{\ac}}}X^{-,+}(E/F_\ac)\big).
\end{multline*}

\er

\end{document}